\documentclass{article}
\usepackage{graphicx} 
\usepackage{amsmath}
\usepackage{booktabs}
\usepackage{mathtools}
\usepackage{amssymb,epsfig,multirow}
\usepackage{amsbsy}
\usepackage{appendix}
\usepackage{graphicx}
\usepackage{subcaption}
\usepackage{color}
\newtheorem{lem}{Lemma}[section]

\newtheorem{cor}{Corollary}[section]
\newtheorem{prop}{Proposition}[section]

\newtheorem{exam}{Example}[section]

\newtheorem{coj}{Conjecture}[section]
\newtheorem{prob}{Problem}[section]

\newtheorem{obser}{Observation}[section]

\newenvironment{proof}{\noindent {\bf Proof. }}{\hfill $\Box$ \newline\par}

\newcommand{\vett}[1]{\left(\begin{array}{ccc}#1\end{array}\right)}

\newcommand{\Real}{\mathbb{R}}
\newcommand{\scalar}[2]{\left<#1,#2\right>}

\newcommand{\C}{\mathcal{C}}

\begin{document}

\title{Is time-optimal speed planning under jerk constraints a convex problem?}

\author{Luca Consolini,\quad  Marco Locatelli\\
Dipartimento di Ingegneria e Architettura\\ Universit\`a di Parma, Italy}

\maketitle

\begin{abstract}                          % Abstract of not more than 200 words.
We consider the speed planning problem for a vehicle moving along an assigned trajectory, under maximum speed, tangential and lateral acceleration, and jerk constraints. The problem is a nonconvex one, where nonconvexity is due to jerk constraints. We propose a convex relaxation, and we present various theoretical properties. In particular, we show that the relaxation is exact under some assumptions. Also, we rewrite the relaxation as a Second Order Cone Programming (SOCP) problem. This has a relevant practical impact, since solvers for SOCP problems are quite efficient and allows solving large instances within tenths of a second.
We performed many numerical tests, and in all of them the relaxation turned out to be exact. For this reason, we conjecture that the convex relaxation is {\em always} exact, although we could not give a formal proof
of this fact.
\end{abstract}

\section{Introduction}

Consider the problem of computing a minimum-time motion of a car-like vehicle from a start configuration
to a target one, while avoiding collisions (obstacle avoidance), and
satisfying kinematic, dynamic, and mechanical constraints (for
instance, on velocities, accelerations and maximal steering
angle).
It is common to solve this problem in two steps. First, we use a geometric path planner to find a suitable path. Then, we perform minimum-time
  speed planning on the planned path (see, for
  instance,~\cite{doi:10.1177/027836498600500304}).
 In this paper, we assume that the path
that joins the initial and final configurations is assigned, and we
aim at finding the time-optimal speed law
that satisfies some kinematic and dynamic constraints.
The problem can be reformulated as an optimization problem, and
it is quite relevant from a practical point of view. In particular, 
in automated warehouses, the speed laws of LGVs are typically planned under acceleration and jerk constraints.
\newline\newline\noindent
In our previous work~\cite{consolini2017scl}, we proposed an optimal
time-complexity
algorithm for finding the time-optimal speed law that satisfies
constraints on maximum velocity and tangential and normal acceleration.
In the subsequent work~\cite{CabConLoc2018coap1}, we included a bound
on the derivative of the acceleration with respect to the arc-length, which results in a convex optimization problem.
Then, in~\cite{9566295}, we considered the
presence of jerk constraints (constraints on the time
derivative of the acceleration). The resulting optimization problem
is non-convex and, for this reason, is significantly more
complex than the ones we discussed in~\cite{consolini2017scl}
and~\cite{CabConLoc2018coap1}.
This work addresses the same problem as~\cite{9566295}. Namely,
we compute a time-optimal control law, taking into account constraints on maximum speed, tangential and lateral acceleration, and jerk. However, we use a completely different approach. Basically, we propose a convex relaxation of the original non-convex problem. Then, we show that, under some assumptions, this relaxation is exact. We reformulate the relaxed problem as a second order cone programming (SOCP) problem. This allows solving the problem very efficiently with modern solvers.

\section{Problem description}
%First, we motivate the study of the following problem, which is the focus of this work.
%\begin{equation}
%\label{eq:orig_int}
%\begin{array}{lll}
%\min_{w} & f(w)=\sum_{i=2}^{n-1} \frac{h}{\sqrt{w_i}} & \\ [6pt]
%& \frac{h}{\sqrt{w_i}}\geq \frac{w_{i-1} - 2 w_i +w_{i+1}}{hJ} & i=2,\ldots,n-1 \\ [6pt]
%& \frac{h}{\sqrt{w_i}}\geq \frac{-w_{i-1} + 2 w_i -w_{i+1}}{hJ} & i=2,\ldots,n-1 \\ [6pt]
%& w_{i+1}-w_i\leq A h & i=2,\ldots,n-1 \\ [6pt]
%& w_{i}-w_{i+1}\leq A h & i=2,\ldots,n-1 \\ [6pt]
%& w_1=w_n=0 &  \\ [6pt]
%& 0\leq w_i\leq w_i^{\max} & i=1,\ldots,n.
%\end{array}
%\end{equation}
%
%Here $w_i$, $i=1,\ldots,n$ is the optimization variable, and $A,J,w_i^{\max}$, $i=1,\ldots,n$ are assigned parameters. We show that Problem~(\ref{eq:orig_int}) arises naturally from the speed planning problem for mobile vehicles.

%\subsection{Speed planning for mobile vehicles}
%
This section is adapted from~\cite{9566295}. For a more detailed discussion, we refer the reader to this reference.
Let $\boldsymbol{\gamma}: [0,s_f] \to \mathbb{R}^2$ be a smooth
function. The image set $\boldsymbol{\gamma}([0,s_f])$ is the path to be
followed, $\boldsymbol{\gamma}(0)$ the initial
configuration, and $\boldsymbol{\gamma}(s_f)$ the final one. % (see Figure \ref{fig:path-to-follow}).
Function $\boldsymbol{\gamma}$ has arc-length parameterization, that is, it is such that
\mbox{$(\forall \lambda \in [0,s_f])$,  $\|\boldsymbol{\gamma}'(\lambda)\| =1$}. 
In this way, $s_f$ is the path length.
%The path is assumed to be parameterized according to its arc
%ength, so that $\lambda$ is the curvilinear abscissa of the path.
We want to compute the speed-law that
minimizes the overall transfer time (i.e., the time needed to go from
$\boldsymbol{\gamma}(0)$ to $\boldsymbol{\gamma}(s_f)$).
To this end, let $\lambda: [0,t_f] \to
[0,s_f]$ be a differentiable monotone strictly increasing function, that
represents the vehicle's arc-length position along the curve as a function of time, and
let
$v: [0,s_f]\to [0,+\infty[$ be such that $(\forall t \in
[0,t_f])\ \dot \lambda(t)=v(\lambda(t))$. In this way, $v(s)$
is the derivative of the vehicle arc-length position, which corresponds to the norm of its velocity vector at position $s$. The position of the vehicle as a function of time is given by ${\bf x}:[0,t_f] \to \Real^2, \ {\bf x}(t)=\boldsymbol{\gamma}(\lambda(t))$.
The velocity
and acceleration are given, respectively, by
$$
\begin{array}{ll}
\dot {\bf x}(t)=\boldsymbol{\gamma}'(\lambda(t)) v(\lambda(t)),\\
\ddot {\bf x}(t)=a_T(t) \boldsymbol{\gamma}'(\lambda(t))+ a_N(t) \boldsymbol{\gamma}'^{\perp} (\lambda(t)),\,
\end{array}
  $$
where $a_T(t)= v'(\lambda(t)) v(\lambda(t))$, $a_N(t)=k(\lambda(t))
v(\lambda(t))^2$ are, respectively, the tangential and normal
components of the acceleration (i.e., the projections of the acceleration
vector $\ddot {\bf x}$ on the tangent and the normal to
the curve). Moreover, $\boldsymbol{\gamma}'^{\perp}(\lambda)=\vett{0&-1\\1&0} \boldsymbol{\gamma}'(\lambda)$ is the normal to
vector $\boldsymbol{\gamma}'(\lambda)$, the tangent of
$\boldsymbol{\gamma}'$ at $\lambda$.
%For details, see for instance Chapter~1.5 of~\cite{do1980differential}.
Here $k:[0,s_f] \to \Real$ is the
scalar curvature, defined as
\mbox{$k(s)=\scalar{\boldsymbol{\gamma}''(s)}{\boldsymbol{\gamma}'(s)^\perp}$}.
Note that $|k(s)|=\|\boldsymbol{\gamma}''(s)\|$.
% Note that the curvature represents
%the velocity with which the tangent vector to the curve rotates with
%respect to the arc length, and quantity $\frac{1}{|k(s)|}$ represents
%the radius of the circle that locally approximates the curve.
In the following, we
assume that $k(s) \in \mathcal{C}^{1}([0,s_f],\Real)$.
%Note that curve $\boldsymbol{\gamma}$ must be chosen such that, for any $s \in [0,
%s_f]$, $|k(s)| \leq R$, where $R$ is the minimum steering radius of
%the vehicle.
The total maneuver time, for a given velocity profile \mbox{$v\in C^1([0,s_f],\Real)$}, is returned by the functional
\begin{equation}
\label{obj_fun_pr}
{\mathcal F}:  C^1([0,s_f],\Real)\rightarrow \Real, \ \ \ {\mathcal F}(v)=\int_0^{s_f} v^{-1}(s) d s.
\end{equation}

We consider the following problem.

\begin{equation}
\label{eqn_problem_pr}
\min_{v \in {\cal V}}  {\cal F}(v),
\end{equation} 
%\begin{subequations}
%\label{eqn_problem_pr}
%\begin{align}
%\min_{v \in C^1([0,s_f],\Real)} & F(v) &\label{obj_fun_pr}\\
%\textrm{subject to, }&v(0)=0,\,v(s_f)=0, \label{inter_con_pr}\\
%& 0< v(s) \leq  \bar v,& s \in ]0,s_f[, \label{con_speed_pr}\\
%& |2 v'(s)v(s)| \leq \theta, &s \in [0,s_f],  \label{con_at_pr}\\
%& |k(s)| v(s)^2 \leq A_N, &s \in [0,s_f],  \label{con_an_pr}
%\end{align}
%\end{subequations}
where the feasible region ${\cal V}\subset  C^1([0,s_f],\Real)$ is defined by the following set of constraints
\begin{subequations}
\label{eqn_problem_constraints}
\begin{align}
v(0)=0,\,v(s_f)=0, \label{inter_con_pr}\\
0\leq v(s) \leq  v_{\max},\ \  s \in ]0,s_f[, \label{con_speed_pr}\\
|v'(s)v(s)| \leq \frac{A}{2}, \ \ s \in [0,s_f],  \label{con_at_pr}\\
  |k(s)| v(s)^2 \leq A_N, \ \ s \in [0,s_f],  \label{con_an_pr}\\
  v''(s) v(s)^2+ v'(s)^2 v(s) \leq \frac{J}{2}, \label{con_jerk_pr}
\end{align}
\end{subequations}
where:
\begin{itemize}
\item  Constraints~\eqref{inter_con_pr} are the initial and final
interpolation
conditions;
\item $k$ is the path curvature;
\item $v_{\max}$, $\frac{A}{2}$, $A_N$, are upper bounds for the velocity, the tangential acceleration, and the normal acceleration, imposed through constraints (\ref{con_speed_pr}), (\ref{con_at_pr}), (\ref{con_an_pr}), respectively;
\item Constraints (\ref{con_jerk_pr}) impose an upper bound $\frac{J}{2}$ on the time derivative of the acceleration (also called ``jerk''). Indeed, note that
\[
\begin{array}{l}
  \frac{d^2}{d t^2} v(\lambda(t))=\frac{d}{d t} (v'(\lambda(t)) \dot \lambda(t))=\frac{d}{d t} (v'(\lambda(t)) v(\lambda(t)))= \\ [8pt]
=v''(\lambda(t)) v(\lambda(t))^2 + v'(\lambda(t))^2 v(\lambda(t))=v''(s) v(s)^2+ v'(s)^2 v(s).
\end{array}
\]
\end{itemize}
After setting $w=v^2$, and noting that $w'=2 v' v$, $w''=2 (v'' v+ v'^2)$, we end up with the following minimum-time problem:

\begin{prob}[Smooth minimum-time velocity planning problem: continuous version]
	\label{cap4:prob:continuos}
	\begin{align}
	&  \min_{w \in C^2} {\displaystyle\int_0^{s_f} w(s)^{-1/2} \, ds} \nonumber\\
%	\intertext{s.t.}
	&w(0)=0, \quad w(s_{f})=0 , \nonumber \\\
	&0  \leq  w(s) \leq \mu^+(s), & s \in [0,s_{f}],  \nonumber\\
	&\frac{1}{2}\left|w^{\prime}(s)\right| \leq \frac{A}{2}, & s\in [0,s_{f}], \nonumber \\ %\label{cap4:acc_bound_cont}\\
	&\frac{1}{2}\left|w^{\prime\prime}(s) \sqrt{w(s)}\right| \le \frac{J}{2},  & s\in[0,s_f],\label{cap4:jerk_bound_cont}
	\end{align}
\end{prob}
where $\mu^{+}$ is the square velocity upper bound, depending on the path curvature, i.e.,
\[
\mu^{+}(s) = \min	\left\{v_{\max}^2,\frac{A_N}{|k(s)|}\right\}.
\] 
Note that the jerk constraint~\eqref{cap4:jerk_bound_cont} is non-convex.
The continuous problem is discretized as follows. We subdivide the path into $n-1$ intervals of equal length (i.e.,
we evaluate function $w$ at points
%$$
%s_i=\frac{(i-1) s_f}{n-1},\ \ \ i=1,\ldots,n,
%$$
$s_i=\frac{(i-1) s_f}{n-1}$, $i=1,\ldots,n$),
so that we have the following $n$-dimensional vector of variables
$$
w=(w_1, w_2,\ldots,w_n)=\left(w(s_1), w(s_2), \ldots, w(s_n)\right).
$$
Then, the finite dimensional version of the problem is: 
\begin{prob} [Smooth minimum-time velocity planning problem: discretized version]
	\label{cap4:prob_disc}
	\begin{align}
	& \qquad \min_{w \in \Real^n} \sum_{i=2}^{n-1} \frac{h}{\sqrt{w_{i}}}
	\label{cap4:obj:disc}\\
%	\intertext{s.t.}
          & w_1=w_n=0 &  \nonumber \\ %\label{cap4:con:initcond}\\ 
	&0  \leq  w_i \leq w_i^{\max} &  i=2, \dots ,n-1, \nonumber \\ % \label{cap4:con:bound}\\ 
	& w_{i+1} - w_i \leq hA, &  i=1, \dots ,n-1,\label{cap4:con:acc}\\ % i=1, \dots ,n-1,
	&w_{i} - w_{i+1} \leq hA, &  i=1, \dots ,n-1,\,\label{cap4:con:dec} \\
	&(w_{i-1} - 2w_i + w_{i+1})\sqrt{w_i}\le h^2J,&  i=2, \dots ,n-1, \label{cap4:con:par}\\
	&-(w_{i-1} - 2w_i + w_{i+1})\sqrt{w_i}\le h^2J,&  i=2, \dots ,n-1, \label{cap4:con:nar}
	\end{align}
\end{prob}
where $w_i^{\max}=\mu^+(s_i)$, for $i=2,\ldots,n-1$. %and, in particular, $u_1 = 0$ and $u_n = 0$, since we are assuming that the initial and final velocity are equal to 0.  
%\textcolor{blue}{In fact, in order to avoid the difficulties related to the nondifferentiability of the square root at 0, a small positive lower bound $\epsilon$ is imposed for each variable $w_i$ and $u_1, u_n$ are imposed to be equal to such value (we set $\epsilon=10^{-7}$ in our experiments).} 
The objective function (\ref{cap4:obj:disc}) is an approximation
of the objective function of Problem \ref{cap4:prob:continuos}, given by a Riemann sum. 
Constraints (\ref{cap4:con:acc}) and (\ref{cap4:con:dec}) are obtained
by a finite difference approximation of $w^{\prime}$.
Constraints~\eqref{cap4:con:par} and~\eqref{cap4:con:nar} 
are obtained by using a second-order central finite
difference to approximate $w^{\prime\prime}$.
Due to jerk constraints~\eqref{cap4:con:par} and~\eqref{cap4:con:nar},
Problem~\ref{cap4:prob_disc} is non-convex. In what follows, we refer to constraints~\eqref{cap4:con:par} as {\em positive jerk constraints}, and to constraints~\eqref{cap4:con:nar} as {\em negative jerk constraints}.  % and
%cannot be solved with the algorithm presented
%in~\cite{CabConLoc2018coap1}. 
After a simple rewriting of the jerk constraints (\ref{cap4:con:par}) and (\ref{cap4:con:nar}), we end up with the following equivalent formulation of Problem~\ref{cap4:prob_disc}:
\begin{equation}
\label{eq:orig_int}
\begin{array}{lll}
\min_{w} & f(w)=\sum_{i=2}^{n-1} \frac{h}{\sqrt{w_i}} & \\ [6pt]
& \frac{h}{\sqrt{w_i}}\geq \frac{w_{i-1} - 2 w_i +w_{i+1}}{hJ} & i=2,\ldots,n-1 \\ [6pt]
& \frac{h}{\sqrt{w_i}}\geq \frac{-w_{i-1} + 2 w_i -w_{i+1}}{hJ} & i=2,\ldots,n-1 \\ [6pt]
& w_{i+1}-w_i\leq A h & i=2,\ldots,n-1 \\ [6pt]
& w_{i}-w_{i+1}\leq A h & i=2,\ldots,n-1 \\ [6pt]
& w_1=w_n=0 &  \\ [6pt]
& 0\leq w_i\leq w_i^{\max} & i=1,\ldots,n.
\end{array}
\end{equation}

\subsection{Main results}
After setting 
\begin{equation}
\label{eq:deltawi}
\Delta w_i=\frac{w_{i-1} - 2 w_i +w_{i+1}}{hJ},
\end{equation}
we have the following equivalent reformulation of problem~(\ref{eq:orig_int}):
\begin{equation}\label{eq:orig_int_equiv}
\begin{array}{lll}
\min_{w,t} & g(t)=\sum_{i=2}^{n-1} t_i& \\ [6pt]
&  t_i\geq \Delta w_i & i=2,\ldots,n-1 \\ [6pt]
&  t_i\geq -\Delta w_i& i=2,\ldots,n-1 \\ [6pt]
& t_i=\frac{h}{\sqrt{w_i}} & i=2,\ldots,n-1 \\ [6pt]
& w_{i+1}-w_i\leq Ah & i=2,\ldots,n-1 \\ [6pt]
& w_{i}-w_{i+1}\leq Ah & i=2,\ldots,n-1 \\ [6pt]
& w_1=w_n=0 &  \\ [6pt]
& 0\leq w_i\leq w_i^{\max} & i=1,\ldots,n.
\end{array}
\end{equation}
It is obtained by adding variables $t_i$, $i=2,\ldots,n-1$, and splitting each pair of jerk constraints into a triple of equivalent constraints. Nonconvexity of this formulation is restricted to the equality constraints $t_i=\frac{h}{\sqrt{w_i}}$ for $i=2,\ldots,n-1$. Then, we can relax the problem into a convex one by replacing the equality constraints with inequalities $t_i\geq \frac{h}{\sqrt{w_i}}$:
\begin{equation}
\label{eq:relax}
\begin{array}{lll}
\min_{w,t} & g(t)=\sum_{i=2}^{n-1} t_i& \\ [6pt]
&  t_i\geq \Delta w_i & i=2,\ldots,n-1 \\ [6pt]
&  t_i\geq -\Delta w_i & i=2,\ldots,n-1 \\ [6pt]
& t_i\geq \frac{h}{\sqrt{w_i}} & i=2,\ldots,n-1 \\ [6pt]
& w_{i+1}-w_i\leq Ah & i=2,\ldots,n-1 \\ [6pt]
& w_{i}-w_{i+1}\leq Ah & i=2,\ldots,n-1 \\ [6pt]
& w_1=w_n=0 &  \\ [6pt]
& 0\leq w_i\leq w_i^{\max} & i=1,\ldots,n.
\end{array}
\end{equation}

Problem~(\ref{eq:relax}), and its relation with~(\ref{eq:orig_int_equiv}), is the main focus of this paper. Despite its simplicity, to our knowledge, relaxation~(\ref{eq:relax}) is new. Our main results are the following ones.
\begin{itemize}
\item We state some properties of relaxation~(\ref{eq:relax}). Namely, we show that a solution of~(\ref{eq:relax}) never violates negative jerk constraints (Proposition~\ref{prop:negjerk}). We show that it can violate positive jerk constraints only if the velocity is equal to its upper bound (Proposition~\ref{prop:posjerk}). Finally, we present a sufficient condition under which the relaxation is exact (Corollary~\ref{cor:1}). 
\item We discuss some properties related to the dual Lagrangian problem of relaxation~(\ref{eq:relax}). In particular, we present a sufficient exactness condition (Proposition~\ref{prop:2}).
  \item We propose a reformulation of relaxation~(\ref{eq:relax}) as a SOCP (see~(\ref{eq:relax_SOCP})). This can be efficiently handled by modern solvers.
\item We present various numerical experiments. In these experiments, relaxation~(\ref{eq:relax}) is \emph{always} exact. This led us to formulate Conjecture~\ref{coj:whiteflag}, is which we surmise that this is always the case. However, \emph{we do not have a proof of this fact}.
  \item We present other numerical experiments, on a generalization of Problem~(\ref{eq:orig_int}) (Section~\ref{sec:var_acc_jerk}), in which acceleration and jerk constraints depend on step $i$ (i.e., the position along the curve). We found that, in many cases, relaxation~(\ref{eq:relax}) is still exact. In all cases, it can be used to find a precise bound on the optimal value of~(\ref{eq:orig_int}). Further, the relaxed solution can be used as a starting point for a local search procedure for non-convex problem~(\ref{eq:orig_int_equiv}).
  \end{itemize}

To our knowledge, all these results are new, since relaxation~(\ref{eq:relax}) and its properties have not been discussed in literature. From a practical point of view, we stress that the proposed non-convex relaxation allows solving very efficiently speed planning problems with jerk constraints. Indeed, our computational times for problems with $1000$ samples are in the order of $0.2$ seconds (see Section~\ref{sec:perftest}).

\subsection{Comparison with existing literature}
For a summary of existing literature, we also refer the reader to our previous paper~\cite{9566295} that, as said, addressed the same problem as this work.
Various works consider problem~(\ref{eq:orig_int}), or similar ones, related to minimum-time speed planning in presence of jerk constraints. Since contraints on maximum jerk are non-convex, these works often use iterative methods to find a local minimum. For instance, in~\cite{debrouwere2013time}, the authors observe that jerk constraints are non-convex, but can be written as the difference of two convex functions. Based on this observation, the authors solve the problem by a sequence of convex subproblems, obtained by linearizing at the current point the concave part of the jerk constraints. In~\cite{Singh15}, the authors reformulate the problem in such a way that its objective function is convex quadratic, while non-convexity lies in difference-of-convex functions. The resulting problem is tackled through the solution of a sequence of convex subproblems obtained by linearizing the concave part of the non-convex constraints. 
Other approaches dealing with jerk constraints do not rely on the solution of convex subproblems. For instance, in~\cite{MacFarlane03}, 
a concatenation of
fifth-order polynomials is employed to provide smooth trajectories, which results in quadratic jerk profiles, while in \cite{Haschke08} cubic polynomials are employed, resulting in piecewise constant jerk profiles.
A very recent and interesting approach to the problem with jerk
constraints is~\cite{pham2017structure}. In this work an approach
based on numerical integration is discussed.
Numerical integration has been first applied under acceleration constraints
in \cite{Bobrow85,Shin85}. In~\cite{pham2017structure} jerk
constraints are taken into account. The algorithm detects a position
$s$ along the trajectory where the jerk constraint is singular, that
is, the jerk term disappears from one of the
constraints. Then, it computes
the speed profile up to $s$ by computing two maximum jerk profiles and then
connecting them by a minimum jerk profile, found by a shooting method. In
general, the overall solution is composed of a sequence of various
maximum and minimum jerk profiles.
This approach does not guarantee reaching a local
minimum of the traversal time.
\newline
Some algorithms use heuristics to quickly find suboptimal solutions
of acceptable quality.
For instance, \cite{Villagra-et-al2012} proposes an algorithm that applies to curves composed of clothoids, circles and
straight lines. The algorithm does not guarantee local
optimality of the solution. Reference~\cite{RaiCGL2019jerk} presents an
efficient heuristic algorithm. Also, this method guaranteea neither global nor local optimality.
Various works in literature consider jerk bounds in the speed optimization problem for
robotic manipulators instead of mobile vehicles. This is a slightly
different problem, but mathematically similar to Problem~\eqref{cap4:prob:continuos}.
In particular, paper~\cite{DONG20071941} presents a method based on the solution of many non-linear and non-convex
subproblems. The resulting algorithm is slow, due to the large number
of subproblems; moreover, the authors do not prove its convergence.
Reference~\cite{ZHANG2012472} proposes a similar method that gives
a continuous-time solution. Again, the method is computationally slow, since
it is based on the numerical solution of many
differential equations; moreover, the paper does not contain a proof of
convergence or of local optimality. In \cite{Palleschi19}, the problem of speed planning for robotic manipulators with jerk constraints is reformulated in such a way that non-convexity lies in simple bilinear terms. Such bilinear terms are replaced by the corresponding convex and concave envelopes, obtaining the so-called McCormick relaxation, which is the tightest possible convex relaxation of the non-convex problem.
Recent works~\cite{9430031} and~\cite{9664638} present efficient heuristics for computing speed profiles with limited jerk.  However, they do not guarantee local or global optimality.

In general, all the above algorithms are able to find good quality solutions, but do not guarantee global (or even local) optimality of the found solution.
In our recent work~\cite{9566295}, we were able to present an algorithm that guarantees local optimality. As said, to our knowledge, the present paper is the only one that guarantees (under some assumptions) global optimality of the found solution.

Some other works replace the jerk constraint
with \emph{pseudo-jerk}, that is the derivative of the acceleration
with respect to arc-length, ending up with a convex optimization problem.
For instance,~\cite{8569414} adds to the objective function a pseudo-jerk penalizing
term. This approach is computationally convenient, but may be overly
restrictive at low speeds. Similarly, works~\cite{9420240,9812155} consider a convex problem obtained by linearizing, or approximating, the jerk constraint. These are convenient approaches from a computational point, but, obviously, do not provide the solution of the original problem.

\section{An alternative way to derive the proposed convex relaxation}
\label{sec:basicidea}
In this section, we show that relaxation~(\ref{eq:relax}) can be obtained from the Lagrangian relaxation of the jerk constraints in~(\ref{eq:orig_int}).
We present the main idea in a slightly more general setting. Consider the following problem.
\begin{prob}
  \begin{equation}
    \label{eqn_prob_gen}
\begin{array}{ll}
  \min_{x \in \C} &\sum_{i\in I}  f_i(x)\\
\textrm{s. t.}& g_i(x)\leq f_i(x), i\in I\,.
\end{array}  
\end{equation}
where for $i\in I$, $f_i,g_i:\Real^n \to \Real$ are convex, and $\C \subset \Real^n$ is convex and compact.
\end{prob}
Note that Problem (\ref{eq:orig_int}) falls into the class of problems (\ref{eqn_prob_gen}) with:
\begin{itemize}
\item $I=\{2,\ldots,n-1\}$ and $x=w$;
\item $\C=\{w\ :\ w_1=w_n=0,\ |w_{i+1}-w_i|\leq A h, \ i=2,\ldots,n-1,\ 0\leq w_i\leq w_i^{\max},\ i=1,\ldots,n\}$;
\item $g_i(w)=\left| \frac{w_{i-1} - 2 w_i +w_{i+1}}{hJ} \right|$, $i=2,\ldots,n-1$.
\end{itemize}
We denote by $F^*$ the optimal value of problem~(\ref{eqn_prob_gen}). Note that this problem is non-convex, since it includes constraints $g_i(x) \leq f_i(x)$. We apply the Lagrangian relaxation to these constraints:
\begin{equation}
    \label{eqn_prob_gen_LR}
\begin{array}{ll}
  \min_{x \in \C} &\sum_{i\in I} f_i(x) + \lambda_i \left(g_i(x)-f_i(x)\right),\\
\end{array}  
\end{equation}
where $\lambda_i$, $i\in I$ are non-negative Lagrange multipliers.
For $\lambda =(\lambda_1,\ldots,\lambda_{|I|}) \in (\Real^+)^{|I|}$, we denote by $F(\lambda)$ be the optimal value of~(\ref{eqn_prob_gen_LR}), and by $x(\lambda)$ a corresponding solution vector.
It is obvious that, for any $\lambda \in (\Real^+)^{|I|}$, $F(\lambda) \leq F^*$, that is, the optimal value of~(\ref{eqn_prob_gen_LR_rew}) is a lower bound for $F^*$. Rewrite~(\ref{eqn_prob_gen_LR}) as
%  \begin{equation}
%    \label{eqn_prob_gen_LR_rew}
$$
\begin{array}{ll}
  \min_{x \in \C} &\sum_{i\in I} (1-\lambda_i) f_i(x) + \lambda_i g_i(x).\\
\end{array}  
$$
%\end{equation}
Note that, if $\lambda_i \in [0,1]$, $i\in I$, this is a convex problem. We consider the following convex relaxation of~(\ref{eqn_prob_gen}):

  \begin{equation}
    \label{eqn_prob_gen_LR_rew}
\begin{array}{ll}
  \max_{\lambda \in [0,1]^{|I|}} \min_{x \in \C} &\sum_{i\in I} (1-\lambda_i) f_i(x) + \lambda_i g_i(x).\\
\end{array}  
\end{equation}

Since $\C$ and $[0,1]^{|I|}$ are compact, we can apply the Von Neumann-Fan minimax theorem, and exchange the maximum and minimum operations, obtaining the following equivalent convex relaxation:

  \begin{equation}
    \label{eqn_prob_gen_LR_rew_exc}
\begin{array}{ll}
  \min_{x \in \C} &\sum_{i\in I} \max\{ f_i(x), g_i(x)\},\\
\end{array}  
\end{equation}
which can also be written as follows:
\begin{equation}
    \label{eqn_prob_gen_LR_rew_t}
\begin{array}{ll}
  \min_{t \in \Real^n,x \in \C} \sum t_i\\
  t_i \geq f_i(x) \\
  t_i \geq g_i(x). \\
  \end{array}  
\end{equation}

Note that relaxation~(\ref{eq:relax}) is exactly the convex relaxation (\ref{eqn_prob_gen_LR_rew_t}) for problem (\ref{eq:orig_int}).

%
%Essentially, bound~(\ref{eqn_prob_gen_LR_rew_t}) is the main convex relaxation used in this work. We can give a sufficient condition under which relaxation~(\ref{eqn_prob_gen_LR_rew}) is exact.
%CONTROLLARE!!! FORSE QUESTO NON LO METTEREI QUI
%
%
%\begin{prop}
% If there exists a solution  $\lambda^*$ of~(\ref{eqn_prob_gen_LR_rew}) such that
%  $\lambda^*_i\leq 1$, $i\in I$.
%then $F(\lambda^*)=F^*$ and $x(\lambda^*)$ is an optimal solution of~(\ref{eqn_prob_gen}).
%\end{prop}
%\begin{proof}
%By contradiction, assume that $F(\lambda^*)< F^*$ or that $x(\lambda^*)$ is not the optimal solution of~(\ref{eqn_prob_gen}).
%Then, $x(\lambda^*)$ is not admissible for~(\ref{eqn_prob_gen}), and there exists $i \in I$ such that
%$g(x_i)>f(x_i)$. This implies that the superderivative of function
%\[
%f(\lambda)= \min_{x \in \C} \sum_{i\in I} (1-\lambda_i) f_i(x) + \lambda_i g_i(x).
%\]
%with respect to $\lambda_i$ is positive, so that $\lambda^*_i=1$.
%\end{proof}.
\section{An exactness condition for the convex relaxation}
Before proceeding, we introduce a slight modification in problem (\ref{eq:orig_int}). 
To simplify the following mathematical analysis, it is worthwhile to modify definition (\ref{eq:deltawi}) as follows:
$$
\Delta w_i=\frac{w_{i-1} - (2+\rho(h)) w_i +w_{i+1}}{hJ},
$$
%\begin{array}{ll}
%\frac{h}{\sqrt{w_i}}\geq \frac{w_{i-1} - (2+\rho(h)) w_i +w_{i+1}}{hJ} & i=2,\ldots,n-1 \\ [6pt]
%\frac{h}{\sqrt{w_i}}\geq \frac{-w_{i-1} + (2+\rho(h)) w_i -w_{i+1}}{hJ} & i=2,\ldots,n-1,
%\end{array}
%$$
where $\rho(h)=o(h^2)>0$. 
Provided that $h$ is small, 
$$
\frac{h}{\sqrt{w_i}}\geq \left| \Delta w_i \right|,\ \ \ i=2,\ldots,n-1,
$$
still represents a correct discretization of the continuos constraints (\ref{cap4:jerk_bound_cont}).
%Recall that
%$$
%\Delta w_i=\frac{w_{i-1} - (2+\rho(h)) w_i +w_{i+1}}{hJ}.
%$$
\newline\newline\noindent
Now, we introduce a condition under which we guarantee that the optimal values of the problems (\ref{eq:orig_int}) and (\ref{eq:relax}) are equal. 
We first introduce the following simple lemma.
\begin{lem}
\label{lem:1}
Any feasible solution $w$ of the original problem (\ref{eq:orig_int}) induces a feasible solution $(w,t)$ with $t_i=\frac{h}{\sqrt{w_i}}$ of the relaxation (\ref{eq:relax}). The two feasible solutions have the same objective function value, i.e., $f(w)=g(t)$.
\end{lem}
As a consequence of Lemma \ref{lem:1}, we have the following primal condition for the exactness of the relaxation.
\begin{obser}
\label{obs:1}
If the optimal solution $(w^\star,t^\star)$ of the convex relaxation  (\ref{eq:relax}) is such that
$w^\star$ is feasible for (\ref{eq:orig_int}), then $w^\star$ is also optimal for (\ref{eq:orig_int}), and the convex relaxation has the same optimal value of the original problem. 
\end{obser}
\begin{proof}
In view of Lemma \ref{lem:1}, the objective function value $g(t^\star)$ of (\ref{eq:relax}) at  $(w^\star,t^\star)$ is equal to
 the objective function value $f(w^\star)$ of (\ref{eq:orig_int}) at $w^\star$. Then, given an optimal solution $\bar{w}$ of  (\ref{eq:orig_int}), we must have
$$
f(w^\star)=g(t^\star)\leq f(\bar{w})\leq f(w^\star),
$$
where the first inequality comes from the fact that the optimal value of the relaxation is a lower bound of the optimal value of (\ref{eq:orig_int}), while the second inequality comes from the fact that $\bar{w}$ and $w^\star$ are an optimal and a feasible solution for problem (\ref{eq:orig_int}), respectively. As a consequence, we must have that all inequalities are equalities, so that $w^\star$ is optimal for (\ref{eq:orig_int}) and problems (\ref{eq:orig_int}) and (\ref{eq:relax})  have the same optimal value. 
\end{proof}
The reverse of Lemma \ref{lem:1} is not true, i.e., given a feasible solution $(\tilde{w},\tilde{t})$ of the convex relaxation (\ref{eq:relax}), $\tilde{w}$ may violate some jerk constraint. E.g., for some $i\in\{2,\ldots,n-1\}$, it may hold that:
\begin{equation}
\label{eq:violation}
\frac{h}{\sqrt{\tilde{w}_i}}< \frac{\tilde{w}_{i-1} - (2+\rho(h)) \tilde{w}_i +\tilde{w}_{i+1}}{hJ}.
\end{equation}
Such opportunity of violating jerk constraints has two conflicting effects on the objective function of the relaxation. On the one hand, violating jerk constraints allows enlarging the feasible set and, thus, to reduce the objective function. More precisely, the projection of the feasible set of (\ref{eq:relax}) over the set of variables $w$, i.e., the set
$$
\left\{w\ :\ (w,t)\ \mbox{is feasible for (\ref{eq:relax})}\right\},
$$
strictly contains the feasible region of (\ref{eq:orig_int}).
On the other hand, a violation has a cost. Indeed, if (\ref{eq:violation}) holds, then we have that
$$
\tilde{t}_i\geq \frac{\tilde{w}_{i-1} - (2+\rho(h)) \tilde{w}_i +\tilde{w}_{i+1}}{hJ}>\frac{h}{\sqrt{\tilde{w}_i}},
$$
i.e., the $i$-th term of the objective function of the relaxation (\ref{eq:relax}) is larger than the corresponding term of the objective function of the original problem (\ref{eq:orig_int}). Therefore, the question is whether the gain obtained from the enlargement of the feasible region is able to counterbalance the cost of the violation. If not, optimal solutions of the relaxation (\ref{eq:relax}) do not violate any jerk constraint and, consequently, are also optimal for (\ref{eq:orig_int}), as stated in Observation \ref{obs:1}.

\subsection{Solution algorithm}
\label{sec:sol_alg}
Observation~\ref{obs:1} motivates the following algorithm for solving~(\ref{eq:orig_int}). First, we solve convex relaxation~(\ref{eq:relax}). To this end, in Section~\ref{sec:socp}, we will present an efficient method based on the reformulation of~(\ref{eq:orig_int_equiv}) as a second-order cone program (SOCP).
If the found optimal solution $(w^\star,t^\star)$ is such that $w^\star$ is feasible, then, by Observation~\ref{obs:1}, $w^\star$ is optimal.
\newline\newline\noindent
In Section~\ref{sec:primalnegandpos}, Corollary~\ref{cor:1}, we will present a sufficient condition for the feasibility (and, hence, optimality) of $w^\star$.
However, as we will discuss in more detail in the following (Section \ref{sec:perftest}), this is always the case in our numerical tests, so that we 
make the following conjecture.
\begin{coj}
\label{coj:whiteflag}
The convex relaxation~(\ref{eq:relax}) is exact, i.e., its optimal value is equal to the optimal value of~(\ref{eq:orig_int}), and given an optimal solution $(w^\star,t^\star)$ for~(\ref{eq:relax}), $w^\star$ is feasible and optimal 
for~(\ref{eq:orig_int}). 
\end{coj}
However, in spite of many attempts to give a formal proof of this conjecture, up to now we have not been able to derive it, apart under the mentioned sufficient condition that we will present in Corollary~\ref{cor:1}. 
%Anyway, we do not have a proof of this fact.
To complete the algorithm, in case $w^\star$ is not feasible, $g(t^\star)$ is a lower bound for the solution of~(\ref{eq:orig_int}). Then, we use $(w^\star,t^\star)$ as the starting condition for a generic nonconvex solver for~(\ref{eq:orig_int_equiv}). Let $(\hat x,\hat w)$ be the obtained (feasible) solution. Then, if $\bar w$ is an optimal solution of~(\ref{eq:orig_int}), we can conclude that 
$f(\bar w) \in [g(t^\star),f(\hat w)]$. We can summarize our solution algorithm as follows.

\begin{enumerate}
\item Set $(w^\star,t^\star)$ as the solution of the convex problem~(\ref{eq:relax}).
  \item If $w^\star$ is feasible for (\ref{eq:orig_int}), then $w^\star$ is the optimal solution of~(\ref{eq:orig_int}).
  \item If $w^\star$ is not feasible, use a generic nonconvex solver for~(\ref{eq:orig_int_equiv}), using initial condition $(w^\star,t^\star)$ and let $(\hat w,\hat t)$ be the obtained solution. Then $\hat w$ is a suboptimal solution, and we can conclude that the true optimal solution $\bar w$ of~(\ref{eq:orig_int}) satisfies $f(\bar w) \in [g(w^\star),f(\hat w)]$.
  \end{enumerate}
As said, in all our numerical tests, $w^\star$ is feasible, so that Step~3 is not necessary. 
%This lead us to the following conjecture.
%\begin{coj}
%\label{coj:whiteflag}
%The convex relaxation~(\ref{eq:relax}) is exact, i.e., its optimal value is equal to the optimal value of~(\ref{eq:orig_int}), and given an optimal solution $(w^\star,t^\star)$ for~(\ref{eq:relax}), $w^\star$ is feasible and optimal 
%for~(\ref{eq:orig_int}). 
%\end{coj}
%However, in spite of many attempts to give a formal proof of this conjecture, up to now we have not been able to derive it, apart under the mentioned sufficient condition that we will present in Corollary~\ref{cor:1}. 
But we cannot remove Step 3 since we do not have a formal proof of Conjecture \ref{coj:whiteflag}. Moreover,
Step 3 is necessary if  we consider a generalization of problem~(\ref{eq:orig_int}), in which constraints $A$, $J$ vary with step $i$ (see Section \ref{sec:perftest}). 
As we will see, in these more general cases there exist instances where $w^\star$ is not feasible for (\ref{eq:orig_int}).
%First, we do not have a theoretical proof that $w^\star$ is feasible (apart under the mentioned sufficient condition that we will present in Corollary~\ref{cor:1}). Second, we will also consider more generic instances of Problem~(\ref{eq:orig_int}), in which constraints $A$, $J$ vary with step $i$ (see METTERE LINK). In these tests, $w^\star$ is not always feasible for (\ref{eq:orig_int}).
  \section{Results on negative and positive jerk constraints}
\label{sec:primalnegandpos}
As previously mentioned, we do not have a formal proof of Conjecture \ref{coj:whiteflag}. However, we can prove some strong theoretical properties for $w^\star$.
\newline\newline\noindent
Let us denote by $X^\star$ the set of optimal solutions of the relaxed problem (\ref{eq:relax}). 
We first prove this result.
\begin{prop}
\label{prop:negjerk}
Given the optimal solution $(w^\star,t^\star)\in X^\star$ of the convex relaxed problem (\ref{eq:relax}), $w^\star$ does not violate any negative jerk constraint.
\end{prop}
\begin{proof} 
Assume, by contradiction, that $(w^\star,t^\star)\in X^\star$ and for some $r\in \{2,\ldots,n-1\}$:
\begin{equation}
\label{eq:viol1}
  t_r^\star=\frac{-w^\star_{r-1} + (2+\rho(h)) w^\star_r -w^\star_{r+1}}{hJ}> \frac{h}{\sqrt{w^\star_r}}.
\end{equation}
Note that we must have $w_r^\star>0$. Indeed, if $w_r^\star=0$, then
$$
-w^\star_{r-1} + (2+\rho(h)) w^\star_r -w^\star_{r+1}>0\ \ \ \Rightarrow \ \ \ w^\star_{r-1} + w^\star_{r+1}<0,
$$
which is not possible.
For $\delta>0$ small enough, let us consider the new feasible solution $\bar{w}$ defined as follows:
$$
\bar{w}_r=w_r^\star-\delta,\ \ \ \bar{w}_i=w_i^\star,\ \ i\neq r.
$$
Obviously, the new solution does not violate the constraints $w_i\leq w_i^{\max}$, and, for $\delta$ small enough, it does not violate the nonnegativity constraints and the positive jerk constraints. It also does not violate acceleration constraints. Indeed, 
by (\ref{eq:viol1}) we have that for $h$ small enough, so that $\frac{h}{\sqrt{w^\star_r}}-\rho(h)w^\star_r>0$, it holds that:
$$
\begin{array}{l}
 Ah\geq w^\star_r -w^\star_{r+1}> w^\star_{r-1} -w^\star_{r} +\frac{h^2}{\sqrt{w^\star_r}}-\rho(h)w^\star_r> w^\star_{r-1} -w^\star_{r}\\ [6pt]
Ah\geq w^\star_r -w^\star_{r-1}> w^\star_{r+1} -w^\star_{r}+\frac{h^2}{\sqrt{w^\star_r}}-\rho(h)w^\star_r> w^\star_{r+1} -w^\star_{r},
\end{array}
$$
so that, for a small enough $\delta$, it also holds that:
$$
\begin{array}{l}
 Ah> \bar{w}_r -\bar{w}_{r+1}> \bar{w}_{r-1} -\bar{w}_{r} \\ [6pt]
Ah> \bar{w}_r -\bar{w}_{r-1}> \bar{w}_{r+1} -\bar{w}_{r},
\end{array}
$$
and the solution $\bar{w}$ fulfills all the acceleration constraints.
We first assume that $r>2$.
We have that
$$
\bar{t}_{r}=t_r^\star-\frac{(2+\rho(h))\delta}{hJ},
$$
while
$$
\begin{array}{l}
\bar{t}_{r-1}=\max\left\{  \frac{w^\star_{r-2} - (2+\rho(h))w^\star_{r-1} +w^\star_{r}-\delta}{hJ} , \frac{-w^\star_{r-2} + (2+\rho(h))w^\star_{r-1} -w^\star_{r}+\delta}{hJ}, \frac{h}{\sqrt{ w^\star_{r-1} }}\right\}
\\ [8pt]
\bar{t}_{r+1}=\max\left\{  \frac{w^\star_{r}-\delta - (2+\rho(h)) w^\star_{r+1} +w^\star_{r+2}}{hJ} , \frac{-w^\star_{r}+\delta + (2+\rho(h)) w^\star_{r+1} -w^\star_{r+2}}{hJ}, \frac{h}{\sqrt{ w^\star_{r+1} }}\right\}.
\end{array}
$$
Now, if 
$$
\frac{-w^\star_{r-2} + (2+\rho(h)) w^\star_{r-1} -w^\star_{r}}{hJ} < t_{r-1}^\star,
$$
we have that, for $\delta$ small enough, 
$$
\bar{t}_{r-1}= t_{r-1}^\star,
$$
Instead, if 
$$
\frac{-w^\star_{r-2} + (2+\rho(h)) w^\star_{r-1} -w^\star_{r}}{hJ} = t_{r-1}^\star,
$$
then 
$$
\bar{t}_{r-1}= t_{r-1}^\star+ \frac{\delta}{hJ},
$$
and
$$
\frac{-\bar{w}_{r-2} + (2+\rho(h)) \bar{w}_{r-1} -\bar{w}_{r}}{hJ}>\frac{h}{\sqrt{\bar{w}_{r-1}}}.
$$
Similar for $\bar{t}_{r+1}$. In all cases we have
$$
\bar{t}_{r}+\bar{t}_{r-1}+\bar{t}_{r+1}\leq t_r^\star+t_{r-1}^\star+t_{r+1}^\star-\frac{\rho(h)\delta}{hJ}<t_r^\star+t_{r-1}^\star+t_{r+1}^\star,
$$
so that optimality is contradicted.
The case $r=2$ can be dealt with in a completely analogous way: since $r=2$, we do not have the updated term
$\bar{t}_{1}$ but only the terms $\bar{t}_{2}$ and $\bar{t}_{3}$. In this case $\bar{t}_2+\bar{t}_{3}<0$, which contradicts optimality.
\end{proof}
This theoretical property is strong. The negative jerk constraint is a nonconvex constraint, but, in spite of that, the optimal solution of the relaxed problem never violates it. If we could prove the same for the positive jerk constraint,
equivalence between (\ref{eq:orig_int}) and its convex relaxation (\ref{eq:relax}) would be established. Unfortunately, we do not have such result. Still we 
can prove a weaker result which restricts the cases where a violation of the positive jerk constraint might occur.
\begin{prop}
\label{prop:posjerk}
Given the optimal solution $(w^\star,t^\star)\in X^\star$ of the convex relaxed problem (\ref{eq:relax}), $w^\star$ can violate the $i$-th positive jerk constraint only if $w_i^\star=w_i^{\max}$.
\end{prop}
\begin{proof} 
We assume, again by contradiction, that $w^\star\in X^\star$ and for some $r\in \{2,\ldots,n-1\}$:
\begin{equation}
\label{eq:viol2}
  \frac{w^\star_{r-1} - (2+\rho(h)) w^\star_r +w^\star_{r+1}}{hJ}> \frac{h}{\sqrt{w^\star_r}},
\end{equation}
and $w^\star_r<w_r^{\max}$.
For $\delta>0$ small enough, let us consider the new feasible solution $\bar{w}$ defined as follows:
$$
\bar{w}_r=w_r^\star+\delta,\ \ \ \bar{w}_i=w_i^\star,\ \ i\neq r.
$$
For $\delta$ small enough the new solution does not violate the constraints $w_i\leq w_i^{\max}$ and the negative jerk constraints, while it obviously does not violate the nonnegativity constraints. It also does not violate acceleration constraints. Indeed, 
by (\ref{eq:viol2}) we have that, for $h$ small enough:
$$
\begin{array}{l}
 Ah\geq w^\star_{r-1} -w^\star_{r}> w^\star_{r} -w^\star_{r+1} +\rho(h) w^\star_r> w^\star_{r} -w^\star_{r+1}\\ [6pt]
Ah\geq w^\star_{r+1} -w^\star_{r}> w^\star_{r} -w^\star_{r-1}+\rho(h) w^\star_r> w^\star_{r} -w^\star_{r-1},
\end{array}
$$
so that, for a small enough $\delta$, it also holds that:
$$
\begin{array}{l}
 Ah> \bar{w}_{r-1} -\bar{w}_{r}> \bar{w}_{r} -\bar{w}_{r+1} \\ [6pt]
Ah> \bar{w}_{r+1} -\bar{w}_{r}> \bar{w}_{r} -\bar{w}_{r-1},
\end{array}
$$
and the solution $\bar{w}$ fulfills all the acceleration constraints.
We only discuss the case $r>2$ (the case $r=2$ can be dealt with in an analogous way).
We have that
$$
\bar{t}_{r}=t_r^\star-\frac{(2+\rho(h))\delta}{hJ},
$$
while
$$
\begin{array}{l}
\bar{t}_{r-1}=\max\left\{  \frac{w^\star_{r-2} - (2+\rho(h)) w^\star_{r-1} +w^\star_{r}-\delta}{hJ} , \frac{-w^\star_{r-2} + (2+\rho(h)) w^\star_{r-1} -w^\star_{r}+\delta}{hJ}, \frac{h}{\sqrt{ w^\star_{r-1} }}\right\}
\\ [8pt]
\bar{t}_{r+1}=\max\left\{  \frac{w^\star_{r}-\delta - (2+\rho(h)) w^\star_{r+1} +w^\star_{r+2}}{hJ} , \frac{-w^\star_{r}+\delta + (2+\rho(h)) w^\star_{r+1} -w^\star_{r+2}}{hJ}, \frac{h}{\sqrt{ w^\star_{r+1} }}\right\}.
\end{array}
$$
Now, if 
$$
\frac{w^\star_{r-2} - (2+\rho(h)) w^\star_{r-1} +w^\star_{r}}{hJ} < t_{r-1}^\star,
$$
we have that, for $\delta$ small enough, 
$$
\bar{t}_{r-1}= t_{r-1}^\star,
$$
Instead, if 
$$
\frac{w^\star_{r-2} -(2+\rho(h)) w^\star_{r-1} +w^\star_{r}}{hJ} = t_{r-1}^\star,
$$
then 
$$
\bar{t}_{r-1}= t_{r-1}^\star+ \frac{\delta}{hJ},
$$
and
$$
\frac{\bar{w}_{r-2} - (2+\rho(h)) \bar{w}_{r-1} +\bar{w}_{r}}{hJ}>\frac{h}{\sqrt{\bar{w}_{r-1}}}.
$$
Similar for $\bar{t}_{r+1}$. In all cases we have
$$
\bar{t}_{r}+\bar{t}_{r-1}+\bar{t}_{r+1}\leq t_r^\star+t_{r-1}^\star+t_{r+1}^\star-\frac{\rho(h)\delta}{hJ}<t_r^\star+t_{r-1}^\star+t_{r+1}^\star,
$$
so that optimality is contradicted.
\end{proof}
We can also prove the following proposition.
\begin{prop}
\label{prop:1}
Given the optimal solution $(w^\star,t^\star)$ of problem (\ref{eq:relax}), we have that $w^\star$ might violate the positive jerk constraint only at some $i$ such that also the maximum speed $w^{\max}$ violates it.
\end{prop}
\begin{proof}
If $w^\star$ violates the $i$-th positive jerk constraint we know that $w_i^\star=w_i^{\max}$, so that
$$
\frac{w_{i-1}^\star-(2+\rho(h))w_i^{\max} +w_{i+1}^\star}{hJ} > \frac{h}{\sqrt{w_i^{\max}}}.
$$
But $w_{i-1}^\star\leq w_{i-1}^{\max}$ and  $w_{i+1}^\star\leq w_{i+1}^{\max}$ imply that
$$
\frac{w_{i-1}^{\max}-(2+\rho(h))w_i^{\max} +w_{i+1}^{\max}}{hJ} > \frac{h}{\sqrt{w_i^{\max}}},
$$
so that also the maximum speed $w^{\max}$ violates the $i$-th positive jerk constraint. 
\end{proof}
Then, we also proved the following corollary.
\begin{cor}
\label{cor:1}
If the maximum speed $w^{\max}$ fulfills all positive jerk constraints, then the relaxation (\ref{eq:relax}) is exact.
\end{cor}
\section{Dual Lagrangian problem}
\label{sec:dual}
In this section we discuss our problems from a dual perspective.
Let:
$$
W=\{w\ :\ 0\leq w_i\leq w_i^{\max},\ i=2,\ldots,n-1,\ w_1=w_n=0\}.
$$
The Lagrangian function for problem (\ref{eq:relax}) is:
$$
\begin{array}{lll}
{\cal L}(w,t,\lambda,\gamma,\xi,\beta,\alpha) & = &  \sum_{i=2}^{n-1} (1-\lambda_i-\gamma_i-\xi_i) t_i +\lambda_i \Delta w_i-\gamma_i\Delta w_i+\xi_i \frac{h}{\sqrt{w_i}} + \\ [6pt]
& + & \alpha_i(w_{i+1}-w_i- Ah)+\beta_i(w_{i}-w_{i+1}- Ah),
\end{array}
$$
and the dual Lagrangian problem is:
\begin{equation}
\label{eq:duallagorig}
\max_{\lambda,\gamma,\xi,\beta,\alpha\geq 0}\ \ \  \min_{w\in W, t} {\cal L}(w,t,\lambda,\gamma,\xi,\beta,\alpha).
\end{equation}
We denote by $W^\star(\lambda,\gamma,\alpha,\beta,\xi)$ the set of optimal solutions of the inner minimization problem.
Since Slater's condition holds for the convex relaxation (\ref{eq:relax}), i.e., its feasible region has a nonempty interior, the Lagrangian dual of the convex relaxation has the same optimal value as the convex relaxation itself.
Since the inner minimization is unbounded from below for $1-\lambda_i-\gamma_i-\xi_i\neq 0$, we can impose the equality or, equivalently, we can replace each $\xi_i$ with $1-\lambda_i-\gamma_i$, and remove variables $t_i$ from the inner problem, so that we can define
$$
\begin{array}{lll}
 {\cal L}'(w,\lambda,\gamma,\beta,\alpha) & = & \sum_{i=2}^{n-1} \frac{h}{\sqrt{w_i}} +
\lambda_i\left(\Delta w_i-\frac{h}{\sqrt{w_i}}\right)+\gamma_i\left(-\Delta w_i-\frac{h}{\sqrt{w_i}}\right)+ \\ [6pt]
  & +& \alpha_i(w_{i+1}-w_i- Ah)+\beta_i(w_{i}-w_{i+1}- Ah),
\end{array}
$$
and the minimization problem reduces to:
\begin{equation}
\label{eq:innermin}
\min_{w\in W}  {\cal L}'(w,\lambda,\gamma,\beta,\alpha),
\end{equation}
whose optimal set is denoted by $W^\star(\lambda,\gamma,\alpha,\beta)$ the set of optimal solutions of the inner minimization problem.
We make the following observation.
\begin{obser}
The dual Lagrangian of the relaxed problem (\ref{eq:relax}) is equivalent to the dual Lagrangian of the original problem (\ref{eq:orig_int}).
\end{obser}
\begin{proof}
It is enough to notice that the Lagrangian function of problem (\ref{eq:orig_int}) is equivalent to ${\cal L}'$.
\end{proof}
We observe that for $1-\lambda_i-\gamma_i<0$, the minimization problem is unbounded from below, so that we further impose that $\lambda_i+\gamma_i\leq 1$. After reorganizing the different terms, the objective function of the minimization problem can be written as the convex separable function:
$$
 \sum_{i=2}^{n-1} \left[(1-\lambda_i-\gamma_i) \frac{h}{\sqrt{w_i}}+(\Delta\lambda_i-\Delta \gamma_i+\beta_i-\alpha_i-\beta_{i-1}+\alpha_{i-1}) w_i -(\beta_i+\alpha_i)hA\right],
$$
where:
$$
\Delta\lambda_i=\frac{\lambda_{i-1}-(2+\rho(h))\lambda_i+\lambda_{i+1}}{hJ},\ \ \ \Delta \gamma_i=\frac{\gamma_{i-1}-(2+\rho(h))\gamma_i+\gamma_{i+1}}{hJ}.
$$
%We can further impose that
%$$
%\Gamma_i(\lambda,\gamma,\alpha,\beta)=\Delta\lambda_i-\Delta \gamma_i+\beta_i-\alpha_i -\beta_{i-1}+\alpha_{i-1}\geq 0,
%$$
%otherwise the problem is unbounded from below.
Now, let us denote by $\omega=(\lambda,\gamma,\alpha,\beta)$ the vector of dual variables
and set
$$
\Gamma_i(\omega)=\Delta\lambda_i-\Delta \gamma_i+\beta_i-\alpha_i -\beta_{i-1}+\alpha_{i-1}.
$$
We consider the following five subsets:
$$
\begin{array}{l}
\Omega^i_1=\left\{\omega\geq 0\ :\ \Gamma_i(\omega)>0,\  \lambda_i+\gamma_i<1,\ \ \mbox{ and } h^{\frac{2}{3}}\left[\frac{1-\lambda_i-\gamma_i}{2\Gamma_i(\omega)}\right]^{\frac{2}{3}} \leq w_i^{\max} \right\} \\ [6pt]
\Omega^i_2=\left\{\omega\geq 0\ :\  \left(\Gamma_i(\omega)\leq 0 \mbox{ and }\lambda_i+\gamma_i<1\right) \mbox{ or } 
\left(\Gamma_i(\omega)>0\mbox{ and } h^{\frac{2}{3}}\left[\frac{1-\lambda_i-\gamma_i}{2\Gamma_i(\omega)}\right]^{\frac{2}{3}} > w_i^{\max}\right)\right\} \\ [6pt]
\Omega^i_3=\left\{\omega\geq 0\ :\  \Gamma_i(\omega)<0 \mbox{ and } \lambda_i+\gamma_i=1\right \} \\ [6pt]
\Omega^i_4=\left\{\omega\geq 0\ :\ \Gamma_i(\omega)>0 \mbox{ and } \lambda_i+\gamma_i=1 \right\} \\ [6pt]
\Omega^i_5=\left\{\omega\geq 0\ :\  \Gamma_i(\omega)=0 \mbox{ and } \lambda_i+\gamma_i=1\right\}.
\end{array}
$$
Then, the solution(s) of the minimization problem, belonging to the optimal set $W^\star(\omega)$, can be given in closed form:
%\begin{equation}
%\label{eq:optsol}
%w_i^\star(\omega)=
%\left\{
%\begin{array}{ll}
%\left[\frac{1-\lambda_i-\gamma_i}{2\Gamma_i(\omega)}\right]^{\frac{2}{3}} & \mbox{if } \Gamma_i(\omega)>0,\  \lambda_i+\gamma_i<1,\ \ \mbox{ and } \left[\frac{1-\lambda_i-\gamma_i}{2\Gamma_i(\omega)}\right]^{\frac{2}{3}} \leq w_i^{\max} \\ [6pt]
%w_i^{\max} &  \mbox{if } \left(\Gamma_i(\omega)=0 \mbox{ and }\lambda_i+\gamma_i<1\right) \mbox{ or } 
%\left(\Gamma_i(\omega)>0\mbox{ and } \left[\frac{1-\lambda_i-\gamma_i}{2\Gamma_i(\omega)}\right]^{\frac{2}{3}} > w_i^{\max}\right)\mbox{ or } \Gamma_i(\omega)<0   \\ [6pt]
%0 & \mbox{if }\Gamma_i(\omega)>0 \mbox{ and } \lambda_i+\gamma_i=1 \\ [6pt]
%[0,w_i^{\max}] & \mbox{if } \Gamma_i(\omega)=0 \mbox{ and } \lambda_i+\gamma_i=1, \\ [6pt]
%\end{array}
%\right.
%\end{equation}
\begin{equation}
\label{eq:optsol}
w_i^\star(\omega)=
\left\{
\begin{array}{ll}
h^{\frac{2}{3}}\left[\frac{1-\lambda_i-\gamma_i}{2\Gamma_i(\omega)}\right]^{\frac{2}{3}} & \mbox{if } \omega\in \Omega^i_1 \\ [6pt]
w_i^{\max} &  \mbox{if } \omega\in \Omega^i_2\cup \Omega^i_3   \\ [6pt]
0 & \mbox{if } \omega\in \Omega^i_4 \\ [6pt]
[0,w_i^{\max}] & \mbox{if } \omega\in \Omega^i_5, \\ [6pt]
\end{array}
\right.
\end{equation}
and its optimal value is:
$$
F(\omega)=\sum_{i=2,\ldots,n-1} F_i(\omega)-(\beta_i+\alpha_i)hA,% \left[\frac{3}{2}\left(1-\lambda_i-\gamma_i\right)^{\frac{2}{3}}\left[2 \Gamma_i(\omega)\right]^{\frac{1}{3}}-(\beta_i+\alpha_i)hA\right].
$$
where
$$
F_i(\omega)=\left\{
\begin{array}{ll}
(1-\lambda_i-\gamma_i) \frac{h}{\sqrt{w_i^{\max}}}+\Gamma_i(\omega) w_i^{\max} & \omega\in \Omega^i_2\cup \Omega^i_3 \\ [6pt]
\frac{3}{2}h^{\frac{2}{3}}\left(1-\lambda_i-\gamma_i\right)^{\frac{2}{3}}\left[2 \Gamma_i(\omega)\right]^{\frac{1}{3}} & \mbox{otherwise}.
\end{array}
\right.
$$
Then, the dual Lagrangian problem is:
\begin{equation}
\label{eq:duallag}
\begin{array}{lll}
\max_{\omega\geq 0} & F(\omega) & \\ [6pt]
        & \lambda_i+\gamma_i\leq 1 & i=2,\ldots, n-1. \\ [6pt]
\end{array}
\end{equation}
Note that $F$ is a continuous and concave function (see, e.g., \cite{Bazar1}).
\subsection{Dual exactness condition}
We observe that from an optimal solution $\omega^\star$ of (\ref{eq:duallag}) we can derive an optimal solution $(\omega^\star,\xi^\star)$ of 
(\ref{eq:duallagorig}) by simply setting $\xi^\star=1-\lambda_i^\star-\gamma_i^\star$ for each $i\in \{2,\ldots,n-1\}$. It is well known (see, e.g., Theorem 6.2.5 in \cite{Bazar1}) that since (\ref{eq:relax}) fulfills the Slater's condition,
it holds that $(\omega^\star,\xi^\star)$ is optimal for (\ref{eq:duallagorig}) and $(w^\star,t^\star)$ is an optimal solution of (\ref{eq:relax}) if and only if:
\begin{itemize}
\item $(\omega^\star,\xi^\star)\geq 0$;
\item $(w^\star,t^\star)$ is feasible for (\ref{eq:relax});
\item $(w^\star,t^\star) \in W^\star(\omega^\star,\xi^\star)$ or, equivalently, $w^\star \in W^\star(\omega^\star)$;
\item all complementarity conditions hold. 
\end{itemize}
Note that, by optimality, we must have
$$
t_i^\star=\max\left\{\frac{h}{\sqrt{w_i^\star}}, \frac{w_{i-1}^\star - (2+\rho(h))w_i^\star +w_{i+1}^\star}{hJ},\frac{-w_{i-1}^\star + (2+\rho(h)) w_i^\star -w_{i+1}^\star}{hJ}\right\}.
$$
Moreover,
if $\xi_i^\star>0$ or, equivalently, $\lambda_i^\star+\gamma_i^\star<1$, the corresponding complementarity condition
$\xi_i^\star\left(t_i^\star-\frac{h}{\sqrt{w_i^\star}}\right)=0$ leads to
$t_i^\star=\frac{h}{\sqrt{w_i^\star}}$, so that $w^\star$ fulfills the jerk constraints at $i$. 
Therefore, a dual exactness condition is the following.
\begin{prop}
If an optimal solution $\omega^\star$ of the dual Lagrangian problem (\ref{eq:duallag}) is such that
$$
\lambda_i^\star+\gamma_i^\star<1,\ \ \ i=2,\ldots,n-1,
$$
then the optimal value of (\ref{eq:duallag}) is equal to the optimal value of (\ref{eq:orig_int}).
\end{prop}
\subsection{Results on negative and positive jerk constraints (dual version)}
Now, let us consider what happens when $\xi_i^\star=0$ for some $i$. In this case
$w^\star$ might violate the $i$-th jerk constraint. 
We first make the following observation.
\begin{obser}
\label{obs:3}
If $\omega^\star\in \Omega_4^i$ for some $i$, then $\omega^\star$ cannot be an optimal solution of (\ref{eq:duallag}).
\end{obser}
\begin{proof}
Assume by contradiction that there exists some optimal solution $\omega^\star$  of (\ref{eq:duallag}) such that $\omega^\star\in \Omega_4^i$ for some $i$. 
Since $\omega^\star\in \Omega_4^i$, we have that
$$
\lambda_i^\star+\gamma_i^\star=1,\ \ \ \Gamma_i(\omega^\star)>0.
$$
Let us consider a further feasible solution $\bar{\omega}$ such that $\bar{\omega}=\omega^\star$ except for
$$
\bar{\lambda}_i=\lambda^\star_i-\delta,
$$
for some $\delta>0$.
We have that
$$
F_i(\omega^\star)=0,\ \ \ F_i(\bar{\omega})\geq \frac{3}{2}h^{\frac{2}{3}}\left[2\Gamma_i(\omega^\star)\right]^{\frac{1}{3}} \delta^{\frac{2}{3}}.
$$
Moreover:
\begin{itemize}
\item if $\omega^\star\in \Omega_4^{i-1}$, then
$F_{i-1}(\bar{\omega})=F_{i-1}(\omega^\star)=0$;
\item if $\omega^\star\in \Omega_2^{i-1}\cup  \Omega_3^{i-1}\cup \Omega_5^{i-1}$, then
$F_{i-1}(\bar{\omega})=F_{i-1}(\omega^\star)-\delta \frac{w_{i-1}^{\max}}{hJ}$;
\item if $\omega^\star\in \Omega_1^{i-1}$, then
$F_{i-1}(\bar{\omega})\approx F_{i-1}(\omega^\star)-\delta \frac{h^{\frac{2}{3}}}{hJ}\left[\frac{1-\lambda^\star_{i-1}-\gamma^\star_{i-1}}{2\Gamma_{i-1}(\omega^\star)}\right]^{\frac{2}{3}}\geq F_{i-1}(\omega^\star)-\delta \frac{w_{i-1}^{\max}}{hJ}$. 
\end{itemize}
Therefore, in every case $F_{i-1}(\bar{\omega})\geq F_{i-1}(\omega^\star)-\delta \frac{w_{i-1}^{\max}}{hJ}$. Similarly,  $F_{i+1}(\bar{\omega})\geq F_{i+1}(\omega^\star)-\delta \frac{w_{i+1}^{\max}}{hJ}$. Now, since for $\delta>0$ and small enough
$$
\begin{array}{lll}
F_{i-1}(\bar{\omega})+F_{i}(\bar{\omega})+F_{i+1}(\bar{\omega}) & \geq &  F_{i-1}(\omega^\star)+F_i(\omega^\star)+F_{i+1}(\omega^\star)+\frac{3}{2}h^{\frac{2}{3}}\left[2\Gamma_i(\omega^\star)\right]^{\frac{1}{3}} \delta^{\frac{2}{3}}-\delta \frac{w_{i-1}^{\max}}{hJ}-\delta \frac{w_{i+1}^{\max}}{hJ} >\\ [8pt]
& > & F_{i-1}(\omega^\star)+F_i(\omega^\star)+F_{i+1}(\omega^\star),
\end{array}
$$
optimality of $\omega^\star$ is contradicted.
\end{proof}
Another implication of the complementarity conditions is that at optimal solutions $(\omega^\star,\xi^\star)$ of (\ref{eq:duallagorig}) it must hold that 
$\lambda_i^\star \gamma_i^\star=0$ for all $i\in \{2,\ldots,n-1\}$. Indeed, if $\lambda_i^\star \gamma_i^\star>0$, we must have $\lambda_i^\star, \gamma_i^\star>0$ and, due to the complementarity conditions:
$$
\begin{array}{l}
t_i^\star=\frac{w_{i-1}^\star - (2+\rho(h)) w_i^\star +w_{i+1}^\star}{hJ} \\ [6pt]
t_i^\star=\frac{-w_{i-1}^\star + (2+\rho(h)) w_i^\star -w_{i+1}^\star}{hJ},
\end{array} 
$$
which is possible only if $w_{i-1}^\star - (2+\rho(h)) w_i^\star +w_{i+1}^\star=0$ and, consequently, $t_i^\star=0$, which is not possible since $t_i^\star\geq \frac{h}{\sqrt{w_i^{\max}}}$.
\newline\newline\noindent
Now, assume that $\lambda_i^\star=1,\ \gamma_i^\star=0$. Then, 
$t_i^\star=\frac{w_{i-1}^\star - (2+\rho(h)) w_i^\star +w_{i+1}^\star}{hJ}>0$, so that:
$$
\begin{array}{l}
hA\geq w_{i-1}^\star -w_{i}^\star>w_{i}^\star-w_{i+1}^\star+\rho(h) w_{i}^\star\ \ \Rightarrow\ \    \beta_i^\star=0 \\ [6pt]
hA\geq w_{i+1}^\star -w_{i}^\star>w_{i}^\star-w_{i-1}^\star+\rho(h) w_{i}^\star\ \ \Rightarrow\ \    \alpha_{i-1}^\star=0, \\ [6pt]
\end{array}
$$
where the implications follow from the complementarity conditions.
Recalling the definition of $\Gamma_i$ we have that
$$
\begin{array}{lll}
\Gamma_i(\omega^\star) & = & \frac{\lambda_{i-1}^\star-(2+\rho(h))\lambda_i^\star+\lambda_{i+1}^\star}{hJ}+\frac{-\gamma_{i-1}^\star+(2+\rho(h))\gamma_i^\star-\gamma_{i+1}^\star}{hJ}+\beta_i^\star-\alpha_i^\star -\beta_{i-1}^\star+\alpha_{i-1}^\star= \\ [6pt]
&=& \frac{\lambda_{i-1}^\star-(2+\rho(h))+\lambda_{i+1}^\star}{hJ}+\frac{-\gamma_{i-1}^\star-\gamma_{i+1}^\star}{hJ}-\alpha_i^\star -\beta_{i-1}^\star\leq \frac{\lambda_{i-1}^\star-(2+\rho(h))+\lambda_{i+1}^\star}{hJ}<0,
\end{array}
$$
where the last inequality follows from $\lambda_{i-1}^\star,\lambda_{i+1}^\star\leq 1$. Therefore, $\omega^\star \in \Omega_3^i$ if $\lambda_i^\star=1$ and, consequently,
$w_i^\star(\omega^\star)=w_i^{\max}$.
\newline\newline\noindent
Similarly, assume that $\lambda_i^\star=0,\ \gamma_i^\star=1$. Then, also recalling that $w_i^\star\leq w_i^{\max}$:
$$
t_i^\star=\frac{-w_{i-1}^\star + (2+\rho(h)) w_i^\star -w_{i+1}^\star}{hJ}\geq \frac{h}{\sqrt{w_i^\star}}\geq \frac{h}{\sqrt{w_i^{\max}}}, 
$$ 
so that:
$$
\begin{array}{l}
hA\geq w_{i}^\star -w_{i+1}^\star\geq w_{i-1}^\star-w_{i}^\star-\rho(h) w_{i}^\star+\frac{h}{\sqrt{w_i^{\max}}} \\ [6pt]
hA\geq w_{i}^\star -w_{i-1}^\star\geq w_{i+1}^\star-w_{i}^\star-\rho(h) w_{i}^\star+\frac{h}{\sqrt{w_i^{\max}}}.
\end{array}
$$
For $h$ small enough we have $-\rho(h) w_{i}^\star+\frac{h}{\sqrt{w_i^{\max}}}>0$, so that
$$
\begin{array}{l}
hA\geq w_{i}^\star -w_{i+1}^\star>w_{i-1}^\star-w_{i}^\star\ \ \Rightarrow\ \    \beta_{i-1}^\star=0\\ [6pt]
hA\geq w_{i}^\star -w_{i-1}^\star>w_{i+1}^\star-w_{i}^\star\ \ \Rightarrow\ \    \alpha_i^\star=0.
\end{array}
$$
Again, recalling the definition of $\Gamma_i$:
$$
\begin{array}{lll}
\Gamma_i(\omega^\star) & = & \frac{\lambda_{i-1}^\star-(2+\rho(h))\lambda_i^\star+\lambda_{i+1}^\star}{hJ}+\frac{-\gamma_{i-1}^\star+(2+\rho(h))\gamma_i^\star-\gamma_{i+1}^\star}{hJ}+\beta_i^\star-\alpha_i^\star -\beta_{i-1}^\star+\alpha_{i-1}^\star = \\ [6pt]
 & = & \frac{\lambda_{i-1}^\star+\lambda_{i+1}^\star}{hJ}+\frac{-\gamma_{i-1}^\star+(2+\rho(h))-\gamma_{i+1}^\star}{hJ}+\beta_i^\star+\alpha_{i-1}^\star\geq \frac{-\gamma_{i-1}^\star+(2+\rho(h))-\gamma_{i+1}^\star}{hJ}>0,
\end{array}
$$
where the last inequality follows from $\gamma_{i-1}^\star,\gamma_{i+1}^\star\leq 1$. Therefore, $\omega^\star \in \Omega_4^i$, and, as stated in Observation \ref{obs:3}, $\omega^\star$ cannot be optimal.
\newline\newline\noindent
Then, we proved the following result.
\begin{prop}
\label{prop:2}
The optimal value of (\ref{eq:relax}) can be strictly lower than the optimal value of (\ref{eq:orig_int}) only if at optimal solutions ${\omega}^\star$ of the dual Lagrangian (\ref{eq:duallag}), it holds that
for some $i\in \{2,\ldots,n-1\}$
$$
\lambda_i^\star=1\ \ \ \mbox{and}\ \ \ \Gamma_i(\omega^\star)< 0,
$$
i.e., ${\omega}^\star\in \Omega^i_3$.
\end{prop}
\section{Counterexamples for more general cases}
Up to now we have not been able to generate any instance for which the optimal value of (\ref{eq:relax}) is strictly lower than the optimal value of (\ref{eq:orig_int}), i.e., we have not been able
to show that Conjecture \ref{coj:whiteflag} is false. However, this is possible for problem classes more general than (\ref{eq:orig_int}).
\subsection{Minimum speed limits}
If, in addition to maximum speed limits, we also add strictly positive minimum speed limits, then Conjecture \ref{coj:whiteflag} is false.  A simple example is the following.
\begin{exam}
Let us set $h,J=1$ and $A=+\infty$. Let us assume that the following lower and upper bounds are imposed for the speeds at $i-1,i,i+1$:
$$
\begin{array}{ll}
w_{i-1}^{\min}=M & w_{i-1}^{\max}=M \\[6pt]
w_{i}^{\min}=0 & w_{i}^{\max}=1 \\[6pt]
w_{i+1}^{\min}=1 & w_{i+1}^{\max}=1. 
\end{array}
$$
The positive jerk constraint
$$
w_{i-1}-(2+\rho(h)) w_i +w_{i+1} \leq \frac{1}{\sqrt{w_i}}
$$
in (\ref{eq:orig_int})
can only be fulfilled, for $M>2$ large enough, when $w_i\leq \bar{w}_i$, where $\bar{w}_i<1$ is such that
$$
\frac{1}{\sqrt{\bar{w}_i}} = M-(2+\rho(h)) \bar{w}_i +1>M-1.
$$
Then, at an optimal solution of (\ref{eq:orig_int}) we have $w_i=\bar{w}_i$ and the contribution of the $i$-term in the objective function 
is $\frac{1}{\sqrt{\bar{w}_i}} $.
Instead, in (\ref{eq:relax}) we can consider a feasible solution with $w_i=1$ and, consequently, the $i$-th term of the objective function is
$$
t_i=M-1+\rho(h)< \frac{1}{\sqrt{\bar{w}_i}}.
$$
\end{exam}
\subsection{Variable acceleration and jerk bounds}
\label{sec:var_acc_jerk}
Up to now we have not been able to prove Conjecture \ref{coj:whiteflag}.
% that the optimal value of the relaxed problem (\ref{eq:relax}) is always equal to the optimal value of the problem (\ref{eq:orig_int}), and that
%an optimal solution of the latter can always be derived from an optimal solution of the former. 
However, in Section \ref{sec:primalnegandpos}
we have been able to prove that optimal solutions of the relaxed problem (\ref{eq:relax}) never violate negative jerk constraints (see Proposition \ref{prop:negjerk}), while violations of positive jerk constraints might occur only at points
where the maximum speed constraint is active (see Proposition \ref{prop:posjerk}). Both propositions are proved by contradiction. More precisely, under the assumption that the optimal solution of (\ref{eq:relax}) violates a jerk constraint, it is shown that the optimal solution can be slightly perturbed in such a way that the perturbed solution is (i) still feasible, and (ii) with lower objective function value, which contradicts optimality.
\newline\newline\noindent
As a next step we might wonder what happens if, in the definition of problem (\ref{eq:orig_int}), we replace the constant acceleration bound $A$ with variable bounds $A_i$, $i=2,\ldots,n-1$.  
In this case, the proofs by contradiction of Propositions \ref{prop:negjerk} and  \ref{prop:posjerk} cannot be applied any more.
Indeed, the perturbed solutions employed in those proofs are not guaranteed to be feasible (some acceleration constraint might be violated).
However, both proofs are still valid if we assume that the acceleration bounds do not vary too quickly. In particular, they are still valid if:
\begin{equation}
\label{eq:condacc}
\left|\frac{A_{i+1}-A_i}{h}\right|\leq \frac{J}{w_i^{\max}},
\end{equation}
i.e., if the variation of the acceleration bound is guaranteed to fulfill the jerk constraints. 
\newline\newline\noindent
Similarly, we might wonder what happens if
we replace the constant jerk bound $J$ with  variable bounds $J_i$, $i=2,\ldots,n-1$. Again, we cannot apply the proofs of Propositions \ref{prop:negjerk} and  \ref{prop:posjerk}. Indeed,
in this case feasibility of the perturbed solution is maintained (provided that the acceleration bound is constant), but the perturbed solution might not have a lower objective function value with respect to the optimal one, thus not leading to a contradiction.
\newline\newline\noindent
As we will see through the experiments in Section \ref{sec:perftest}, in these more general cases we could detect instances where the optimal value of (\ref{eq:relax}) is a strict lower bound of the optimal value of (\ref{eq:orig_int}), except for the case with fixed jerk bound and variable acceleration bounds fulfilling condition (\ref{eq:condacc}). 
Thus, interestingly, the cases for which we have been able to find instances for which the lower bound is strict are also those to which we cannot extend the proofs of Propositions \ref{prop:negjerk} and  \ref{prop:posjerk}.

\section{Numerical tests}
\label{sec:num_tests}
\subsection{SOCP reformulation}
\label{sec:socp}
In the following numerical tests, we will apply the algorithm presented in Section~\ref{sec:sol_alg} to various cases. For efficiency of computation it is convenient to reformulate the relaxed problem~(\ref{eq:relax}) as a SOCP (Second-Order Cone Programming), for which solvers more efficient than generic nonlinear solvers are available.
To this end, note that constraint $t\geq \frac{h}{\sqrt{w}}$ is equivalent to $t^2 w \geq h^2$. Since $t,w \geq 0$, following~\cite{alizadehSecondorderConeProgramming2003}, this last constraint is equivalent to
\[
  \begin{array}{ll}
    x_2^2 \leq t h\\
    x_1^2 \leq t w h\\
    h^2 \leq x_1 x_2\\
    x_1, x_2\geq 0.
    \end{array}
\]
The quadratic constraints can be refomulated as the following SOCP constraints:
\[
  \left\| \vett{ \frac{2x_2}{\sqrt{h}}\\t-1} \right \| \leq t+ 1,\;
  \left\| \vett{ \frac{2x_1}{\sqrt{h}}\\t-w} \right \| \leq t+ w,\;
  \left\| \vett{2h\\x_2-x_1} \right \| \leq x_2+ x_1. \\    
\]

This leads to the following SOCP reformulation of~(\ref{eq:relax}).

\begin{equation}
\label{eq:relax_SOCP}
\begin{array}{lll}
\min_{w,t} & g(t)=\sum_{i=2}^{n-1} t_i& \\ [6pt]
&  t_i\geq -\frac{w_{i-1}-2 w_i + w_{i+1}}{h J_i} & i=2,\ldots,n-1 \\ [6pt]
           &  t_i\geq \frac{w_{i-1}-2 w_i + w_{i+1}}{h J_i} & i=2,\ldots,n-1 \\ [6pt]
&      x_{2,i}^2 \leq t_i h & i=2,\ldots,n-1\\
&    x_{1,i}^2 \leq t_i w_i h  & i=2,\ldots,n-1\\
&    h^2 \leq x_{1,i} x_{2,i} &i=2,\ldots,n-1 \\
&    x_{1,i}, x_{2,i}\geq 0  &i=2,\ldots,n-1\\
& w_{i+1}-w_i\leq A_i h & i=2,\ldots,n-1 \\ [6pt]
& w_{i}-w_{i+1}\leq A_i h & i=2,\ldots,n-1 \\ [6pt]
& w_1=w_n=0 &  \\ [6pt]
& 0\leq w_i\leq w_i^{\max} & i=1,\ldots,n.
\end{array}
\end{equation}

We can solve~(\ref{eq:relax_SOCP}) with efficient commercial solvers, such as MOSEK or GUROBI. Note that we are considering the generic case in which acceleration and jerk bounds depend on $i$.
\subsection{Performed tests}
\label{sec:perftest}
We performed various numerical tests from randomly generated data.
First, we comment on the parameters choice in~(\ref{eq:relax_SOCP}).
Define vectors $J=(J_2,\ldots,J_{n-1}) \in \Real^{n-2} , A=(A_2,\ldots,A_{n-1}) \in \Real^{n-2}$, and, similarly,  $w^{\max},w,t,x_1,x_2 \in \Real^n$. We represent all parameters in~(\ref{eq:relax_SOCP}) with set $(h,A,J,w^{\max})$.
Define
\[\mathcal{S}(h,A,J,w^{\max})=\{(w,t,x_1,x_2) : (w,t,x_1,x_2) \textrm{ is a solution of~(\ref{eq:relax_SOCP}) with parameters } h,A,J,w^{\max}\}.
\]
Then, by substituition, we can verify the following scaling property
\[
  (w,t,x_1,x_2) \in \mathcal{S}(h,A,J,w^{\max})
  \Leftrightarrow (\forall \rho \neq 0) \, (w,\rho t,\rho x_1,\rho x_2) \in \mathcal{S}(\rho h, \rho^{-1}A, \rho^{-2} J,w^{\max}).
\]

Because of this property, it is not restrictive to assume $h=1$, since we can always reduce to this case by scaling the parameters with $\rho=h^{-1}$.

In our tests, we used $h=1$ and $n=1000$. We randomly generated $w^{\max}$ with the following possible strategies:
\begin{itemize}
\item A random vector, in which each component is uniformly distributed in interval $[0.01,100]$ (named \texttt{rnd} in the results table~\ref{tab:results}).
  \item A piecewise constant vector, in which values $w^{\max}$ are constant for each $n/10$ consecutive components. Again, the values of $w^{\max}$ are random numbers, uniformly distributed in $[0.01,100]$ (\texttt{pw cnst}).
\item A piecewise linear function, obtained by linearly interpolating random values at each $n/10$ samples, uniformly distributed in interval $[0.1,100]$ (\texttt{pw lin}).
\end{itemize}

Then, we generated $A$ according to one of the three strategies:
\begin{itemize}
\item A random constant value, obtained from a uniform distribution in interval $[0.1,100]$ (\texttt{cnst}).
\item A random vector, in which each component is uniformly distributed in interval $[0.1,100]$ (\texttt{rnd}).
\item A regularized random vector, fulfilling condition~(\ref{eq:condacc}) (\texttt{reg}).
\end{itemize}

Similarly, we generated $J$ in two possible ways:
\begin{itemize}
\item A random constant value, obtained from a uniform distribution in interval $[0.01,100]$ (\texttt{cnst}).
\item A random vector, in which each component is uniformly distributed in interval $[0.01,100]$ (\texttt{rnd}).
\end{itemize}

 \begin{figure}
	\centering
	\begin{subfigure}[b] {0.3\linewidth}
		\includegraphics[width=\linewidth]{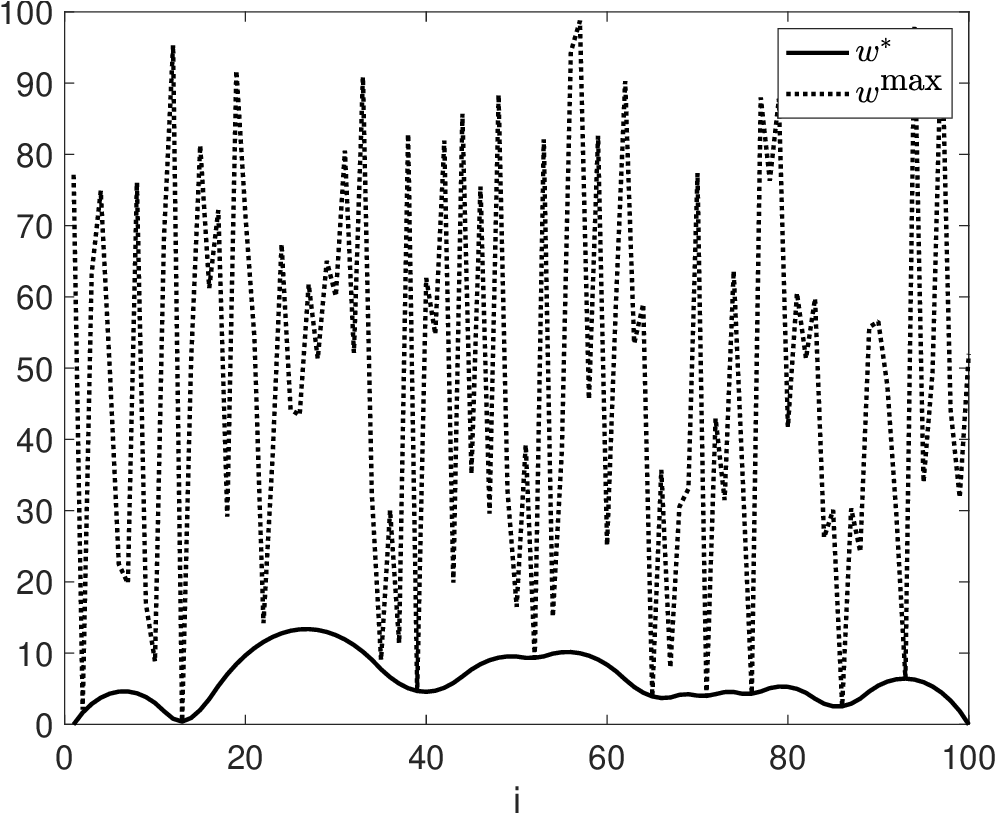}
		\caption{Random $w^{\max}$}
		\label{fig:subfigA}
	\end{subfigure}
	\begin{subfigure}[b]{0.3\linewidth}
		\includegraphics[width=\linewidth]{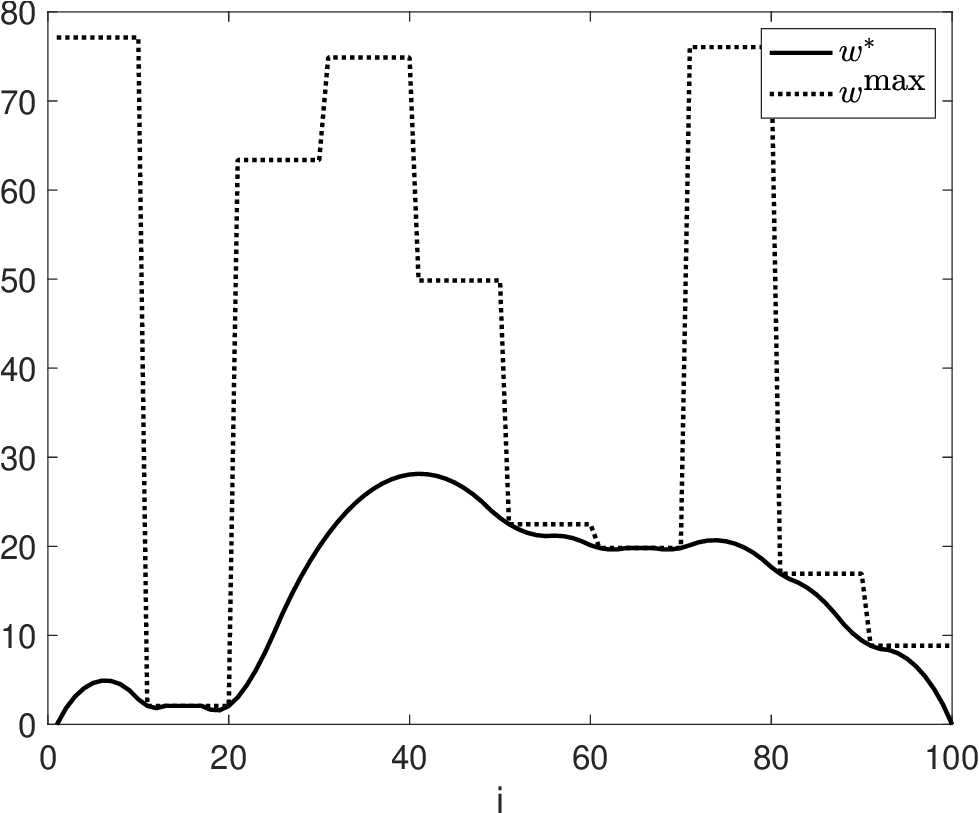}
		\caption{Piecewise constant $w^{\max}$}
		\label{fig:subfigB}
	\end{subfigure}
	\begin{subfigure}[b]{0.3\linewidth}
	        \includegraphics[width=\linewidth]{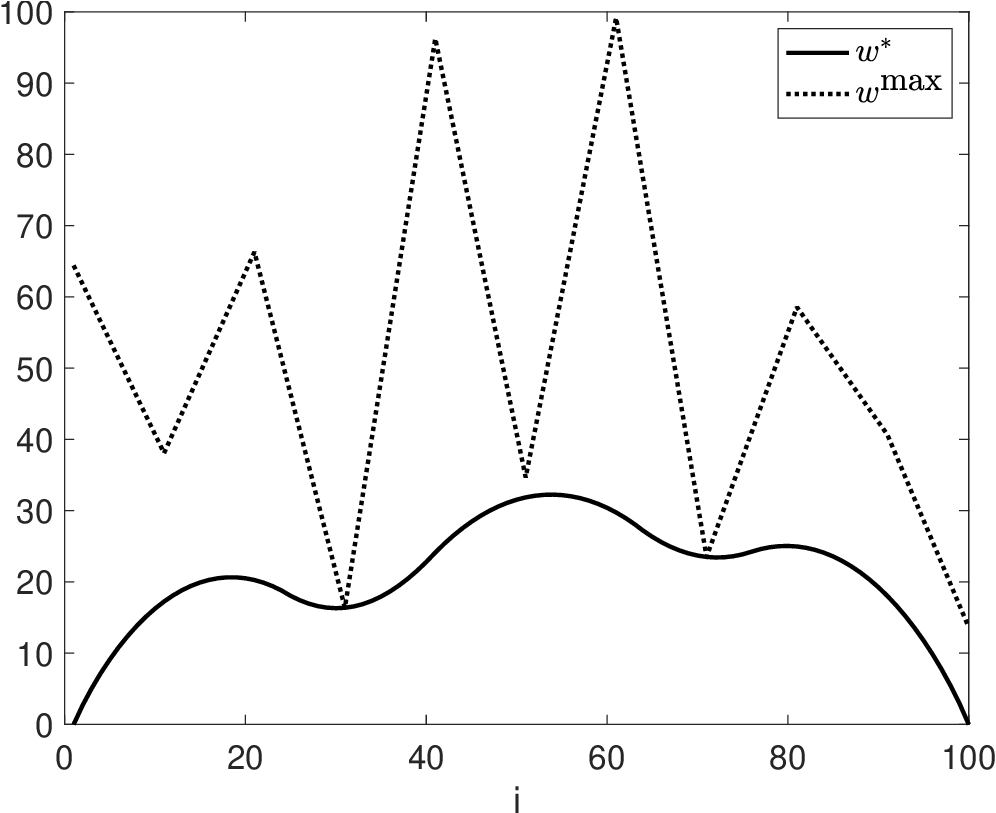}
	        \caption{Piecewise linear $w^{\max}$}
	        \label{fig:subfigC}
         \end{subfigure}
	\caption{Three random problems with $w^{\max}$ generated by different strategies. Parameters $J$ and $A$ are constant, and $N=100$. The dashed line represents $w^{\max}$, the solid one the optimal solution $w^*$.}
	\label{fig:profiles_w}
\end{figure}

Figure~\ref{fig:profiles_w} represents three random problems with $n=100$. We generated $w^{\max}$ with the three proposed strategies (random, piecewise constant, and piecewise linear), while $A$ and $J$ are constant. The dashed lines represent $w^{\max}$, while the solid line is the optimal solution $w^*$.

Given all possible choices, we have a total of $3 \times 3 \times 2=18$ test types. For each test type, we considered $1000$ instances.
We solved each instance with the algorithm presented in Section~\ref{sec:sol_alg}. For the convex relaxation, we used solver MOSEK on the SOCP reformulation~(\ref{eq:relax_SOCP}). We considered a convex relaxation solution $w^*$ a feasible one if
\begin{equation}
  \label{eq:cond_rel}
\max_{i=2,\ldots,N} \left|w^*_{i-1}-2 w^*_i+ w^*_{i+1}\right|-\frac{J_i}{\sqrt{w^*_i}}\leq 10^{-5},
\end{equation}
that is, if the violation of the jerk constraint does not exceed $10^{-5}$.
For those instances in which the convex relaxation solution $w^*$ was not feasible, we used IPOPT to solve problem (\ref{eq:orig_int}) with starting point $w^*$. Calling $\hat w$ the solution obtained by IPOPT, we then computed the relative gap $\frac{f(\hat w)-g(w^*)}{g(w^*)}$. Recall that the (true) optimal solution $\bar w$ of the non-relaxed problem satifies $f(\bar w) \in [g(w^*),f(\hat w)]$. Hence, the relative gap is a measure of the remaining uncertainty on the true optimal solution.

Table~\ref{tab:results} reports the obtained results. In particular,
rows $W$, $A$, $J$ represent the method used for generating the random $w^{\max}$, $A$, $J$ parameter vectors. Then, for each case, we report the number of non-exact relaxations, that is the number of times (over the 1000 tests) in which the convex relaxation violated condition~(\ref{eq:cond_rel}), and the maximum and mean values of the jerk error. Then, for those cases in which we have at least one non-exact relaxation, we report the maximum and mean gap. Finally, we report the MOSEK computation time for the convex relaxation. In these tests, we used a laptop with an 8 cores Apple M1 Pro processor and 16GB of RAM. Note that the computational times are comparable to the ones reported in our paper~\cite{9566295}, where we used a completely different approach that only guarantees local optimality.  
\newline\newline\noindent
We stress the fact that, in all the 3,000 instances with constant value for $A$ and $J$, no failure has been observed. In fact, in an attempt to find a counterexample showing that the optimal values of (\ref{eq:orig_int}) and (\ref{eq:relax}) are not always equal, we performed many more tests with respect to the 3,000 reported in Table~\ref{tab:results}, but no failure has ever been observed. 
No failure has ever been observed also in the case of constant $J$ value and acceleration bounds fulfilling condition~(\ref{eq:condacc}). Very few faliures are observed when $J$ is constant while $A$ is a randomly generated vector.
Instead, the number of failures becomes significantly large when $J$ is a randomly generated vector, especially when also the maximum speed profile is generated in a fully random way. But even in cases where many failures occur, the percentage gap is never too large and never exceeds 1\%.

\section{Conclusions and future research}
In this paper we addressed the speed planning problem along a given trajectory under maximum speed, tangential and lateral acceleration, and jerk constraints. The problem is refomulated as a nonconvex optimization problem, where nonconvexity is due to the jerk constraints. We derived a convex relaxation of such nonconvex problem, exploiting the special structure of the jerk constraints, which can be rewritten in such a way that their left-hand side is equal to the terms of the objective function. Different strong theortical properties of the convex relaxation are proved, in particular, the fact that its optimal solution never violates negative jerk constraints and can only violate positive jerk constraints at points where the maximum speed constraint is active, and the fact that the relaxation is exact if the maximum speed profile does not violate the positive jerk constraints. In fact, we conjectured that the convex relaxation is always exact, i.e., its optimal value is always equal to the optimal value of the speed planning problem, and an optimal solution of the latter can always be derived from an optimal solution of the former. We have been unable to prove this result, but after performing thousands of tests we have never found a counterexample to our conjecture. Instead, the conjecture turns out to be false as soon as we extend the class of problems, in particular, by allowing different bounds for the jerk along the trajectory. It has also been shown that the convex relaxation can be rewritten as a SOCP problem. This has a relevant practical impact, since solvers for SOCP problems are quite efficient and allows solving large instances within tenths of a second.
\newline\newline\noindent
As a possible topic for future research, we are interested in evaluating the performance of the proposed approach to robotic manipulators. These problems can be reformulated in a way similar to
the problems addressed in this paper. We do not expect to be able to extend Conjecture \ref{coj:whiteflag} to such problems but we do expect that, even in cases where the convex relaxation is not exact, the final solution returned by the proposed solution algorithm has a small percentage gap with respect to the optimal one. 

\begin{table}
  \begin{tabular}{lllllllll}

    W & A & J & Non-exact & Max jerk err & Mean jerk err & Max gap[\%] & Mean gap[\%] & Mean time [s] \\ 
\hline 
rnd & cnst & cnst & 0 & 8.5519e-08 & 8.5519e-11 & N/A & N/A & 0.17068 \\ 
pw cnst & cnst & cnst & 0 & 0 & 0 & N/A & N/A & 0.1533 \\ 
pw lin & cnst & cnst & 0 & 0 & 0 & N/A & N/A & 0.20595 \\ 
rnd & rnd & cnst & 0 & 0 & 0 & N/A & N/A & 0.17975 \\ 
pw cnst & rnd & cnst & 4 & 0.16594 & 0.0002596 & 0.004401 & 6.803e-06 & 0.1747 \\ 
pw lin & rnd & cnst & 0 & 7.9644e-08 & 9.7622e-11 & N/A & N/A & 0.1673 \\ 
rnd & reg & cnst & 0 & 0 & 0 & N/A & N/A & 0.19426 \\ 
pw cnst & reg & cnst & 0 & 0 & 0 & N/A & N/A & 0.18469 \\ 
pw lin & reg & cnst & 0 & 0 & 0 & N/A & N/A & 0.17357 \\ 
rnd & cnst & rnd & 756 & 1.0218 & 0.2061 & 0.17484 & 0.0098672 & 0.18201 \\ 
pw cnst & cnst & rnd & 45 & 0.59223 & 0.0080017 & 0.067108 & 0.00034977 & 0.15341 \\ 
pw lin & cnst & rnd & 0 & 0 & 0 & N/A & N/A & 0.15199 \\ 
rnd & rnd & rnd & 788 & 0.96156 & 0.20684 & 0.19155 & 0.010135 & 0.17649 \\ 
pw cnst & rnd & rnd & 251 & 2.0381 & 0.11864 & 0.64832 & 0.011838 & 0.21448 \\ 
pw lin & rnd & rnd & 191 & 2.3979 & 0.097524 & 0.80218 & 0.013251 & 0.15705 \\ 
rnd & reg & rnd & 808 & 1.1569 & 0.22148 & 0.27889 & 0.010681 & 0.19329 \\ 
pw cnst & reg & rnd & 42 & 0.70596 & 0.0080485 & 0.065445 & 0.00038971 & 0.19225 \\ 
pw lin & reg & rnd & 0 & 0 & 0 & N/A & N/A & 0.21015 \\ 
\hline 
\end{tabular}
\caption{Computational results}
  \label{tab:results}

\end{table}

\bibliographystyle{plain}
\bibliography{biblio,VelPlan,VelPlan2}

\end{document}